\documentclass[12pt,twoside,a4paper]{article}
\usepackage{srcltx}
\usepackage{amsmath,amsfonts,amscd,amssymb,latexsym}
\usepackage{xcolor}
\usepackage[latin1]{inputenc}
\usepackage{hyperref}

\newtheorem{theorem}{Theorem}[section]
\newtheorem{lemma}{Lemma}[section]
\newtheorem{corollary}{Corollary}[section]
\newtheorem{remark}{Remark}[section]
\newcommand{\oti}[2]{\mathop \otimes\limits_{#2}^{#1} v}

\newcommand{\ox}[3]{\mathop {#1}\limits_{#3}^{#2} }

\newcommand{\0}{_{\bar 0}}

\newcommand{\1}{_{\bar 1}}
\newcommand{\F}{{\cal F}}

\def\proof{{\noindent\it Proof. }}
\def\qed{\hfill \rule{2.25mm}{2.25mm}\vspace{0,2cm}}
\renewcommand{\=}{\doteq}
\def\qed{\hfill \rule{2.25mm}{2.25mm}\vspace{10pt}}
\renewcommand{\L}{{\cal L}}

\begin{document}

\noindent
\vspace{1 in}
\begin{center}
\noindent {\Large {\bf ON SIMPLE FILIPPOV
SUPERALGEBRAS\\[2mm] OF TYPE $A(0,n)$}}
\vspace{.5cm}

P. D. Beites and A. P. Pozhidaev
\end{center}
\vspace{0,5cm}

\parbox[c]{13cm}{{\large Abstract}$:$
It is proved that there exist no simple finite-dimensional
Filippov superalgebras of type $A(0,n)$ over an
algebraically closed field of characteristic $0$.\\

}

\section{Introduction}

The concept of $n$-Lie superalgebra was
presented by Daletskii and Kushnirevich, in \cite{Dal}, as a natural generalization of the
$n$-Lie algebra notion introduced by Filippov in $1985$ (see \cite{Fil}). Following \cite{GM} and \cite{EnL}, in this article, we use the terms Filippov superalgebra and Filippov algebra instead of $n$-Lie superalgebra and  $n$-Lie algebra, respectively. Filippov algebras were also known before under the names of Nambu Lie algebras and Nambu algebras. As pointed out in \cite{malc} and \cite{genmalc}, Filippov algebras are a particular case of $n$-ary Malcev algebras (generalizing the fact that every Lie algebra is a Malcev algebra). We may also remark that a $2$-Lie superalgebra is simply known as a Lie superalgebra. The description of the finite-dimensional simple Lie superalgebras over an algebraically closed field of characteristic zero was given
 by Kac in \cite{Kac}.

This work is one more step on the way to the classification of
  finite-dimensional simple Fi\-lippov superalgebras over an algebraically
  closed field of characteristic $0$.
  In \cite{LA}, finite-dimensional commutative $n$-ary Leibniz
  algebras over a field of characteristic $0$ were studied by the second author. He showed
  that there exist no simple ones.
  The finite-dimensional simple Filippov algebras over an
  algebraically closed field of characteristic $0$ were classified by Ling in
  \cite{Ling}. Notice that an $n$-ary commutative Leibniz algebra is
  exactly a Filippov superalgebra with trivial even part, and a Filippov
  algebra is exactly a Filippov superalgebra with trivial odd part.
Bearing in mind these facts, we consider the
 $n$-ary Filippov superalgebras
  with $n\geq 3$, and with nonzero even and odd parts.

Let $G$ be a Lie superalgebra. We say that a Filippov superalgebra
$\F$ has {\it type} $G$ if $Inder(\F)\cong G$ (see definitions below).
A description of simple Filippov superalgebras of type $B(m,n)$ was already
obtained in \cite{Poj8}, \cite{PS} and \cite{P7}.
 The same problem concerning Filippov superalgebras
of type $A(m,n)$ with $m=n$ has recently been solved in \cite{BPoj}. The present work represents one more step towards
the classification of finite-dimensional simple Fi\-lippov superalgebras of type
$A(m,n)$ over an algebraically closed field of cha\-racteristic zero. Concretely, we
establish a negative answer to the existence problem of the mentioned superalgebras when $m=0$.

We start recalling some definitions.

 An {\it $\Omega$-algebra} over a field $k$ is a linear
  space over $k$ equipped with a system of multilinear
  algebraic operations $\Omega = \{\omega_i: | \omega_i | = n_i \in
  {\mathbb N}, \ i \in I \}$, where $| \omega_i |$
  denotes the arity of $\omega_i$.

  An {\it $n$-ary Leibniz algebra} over a field $k$ is
  an $\Omega$-algebra $L$ over $k$ with one
  $n$-ary operation $(\cdot, \cdots, \cdot)$ satisfying the identity
 \begin{center}
  $((x_1, \ldots, x_n), y_2, \ldots, y_n) =
  \displaystyle\sum_{i = 1}^{n} (x_1, \ldots, (x_i, y_2, \ldots, y_n), \ldots, x_n).$ \end{center}  If this operation is anticommutative, we obtain the definition of
  \textit{Filippov} (\textit{$n$-Lie}) algebra over a field.

  An  {\it $n$-ary superalgebra} over a field $k$ is
  a $\mathbb{Z}_2$-graded $n$-ary algebra $L=L\0\oplus L\1$ over $k$,
  that is,

\begin{center} if $x_i\in L_{\alpha_i}, \alpha_i\in \mathbb{Z}_2$, then
  $(x_1, \ldots, x_n)\in L_{\alpha_1+\ldots+\alpha_n}$. \end{center}
An {\it $n$-ary Filippov superalgebra} over $k$ is
  an $n$-ary superalgebra $\F=\F\0\oplus\F\1$ over $k$ with one
  $n$-ary operation $[\cdot, \cdots, \cdot]$ satisfying
 \begin{eqnarray}\label{AC}
  &&\!\!\!\! \!\!\!\![x_1, \ldots, x_{i-1},x_i, \ldots, x_n] =
  -(-1)^{p(x_{i-1})p(x_i)}[x_1, \ldots, x_{i},x_{i-1}, \ldots, x_n],\\
  \label{YaI}
&&\!\!\!\!\!\!\!\![[x_1, \ldots, x_n], y_2, \ldots, y_n]  =
  \sum_{i = 1}^{n}(-1)^{p\bar q_i}
   [x_1, \ldots, [x_i, y_2, \ldots, y_n], \ldots, x_n],
 \end{eqnarray}
 where $p(x)=l$ means that $x\in \F_{\bar l}$, $p=\sum_{i=2}^{n}p(y_i),\
  \bar q_i=\sum_{j=i+1}^{n}p(x_j),\ \bar{q}_n=0$. The identities (\ref{AC})
  and (\ref{YaI})
 are called
 the anticommutativity and the generalized Jacobi identity, respectively.
 By (\ref{AC}), we can rewrite (\ref{YaI}) as
\begin{eqnarray}\label{LYaI}
&&\!\!\!\!\!\!\!\![y_2, \ldots, y_n,[x_1, \ldots, x_n]] =
  \sum_{i = 1}^{n}(-1)^{p q_i}
   [x_1, \ldots, [y_2,\ldots, y_n,x_i], \ldots, x_n],
 \end{eqnarray}
where $q_i=\sum_{j=1}^{i-1}p(x_j),\ q_1=0$. Sometimes, instead of using the long term
 ``$n$-ary superalgebra'', we simply say for short ``superalgebra''.
 If we denote by $L_x=L_{(x_1,\ldots,x_{n-1})}$
 the operator of left multiplication
 $L_xy=[x_1,\ldots,x_{n-1},y]$, then, by (\ref{LYaI}), we get
\begin{eqnarray*}\label{OLYaI}
&&[L_y,L_x] =
  \sum_{i = 1}^{n-1}(-1)^{p q_i}
   L(x_1, \ldots,L_yx_i, \ldots, x_{n-1}),
 \end{eqnarray*}
where $L_y$ is an operator of left multiplication and $p$ its parity.
(Here and afterwards, we denote the supercommutator by $[\,,]$).

Let $L=L\0\oplus L\1$ be an $n$-ary anticommutative superalgebra.
A \textit{subsuperalgebra} $B=B\0\oplus B\1$ of the superalgebra $L$,
$B_{\bar i}\subseteq L_{\bar i}$,
is a $\mathbb{Z}_2$-graded vector subspace of $L$ such that $[B, \ldots, B] \subseteq B$.
A subsuperalgebra $I$ of $L$ is called an \it ideal\/ \rm if
$[I,L, \ldots ,L] \subseteq I$.
The subalgebra (in fact, an ideal) $L^{(1)}=[L, \ldots ,L]$ of $L$
is called the \it derived subsuperalgebra\/ \rm of $L$.
Put $L^{(i)}=[L^{(i-1)}, \ldots ,L^{(i-1)}]$, $i\in {\mathbb N}, i>1$.
The superalgebra $L$ is called {\it solvable} if $L^{(k)}=0$ for some $k$.
Denote by $R(L)$ the maximal solvable ideal of $L$ (if it exists). If $R(L)=0$,
 the superalgebra $L$ is called {\it semisimple}.
The superalgebra $L$ is called \it simple\/ \rm if
$L^{(1)} \neq 0$ and $L$ lacks ideals other than $0$ or $L$.

 The article is organized as follows.

In the second section we recall how to reduce the classification problem
of simple Filippov superalgebras to some question about Lie superalgebras,
using the same ideas as in \cite{Ling}. Concretely, we consider an existence problem for some skewsymmetric homomorphisms of semisimple
Lie superalgebras and their faithful irreducible modules. This section is followed with the third one where we collect some definitions and results on Lie superalgebras that we will apply in the two last sections. We also fix some notations with the same purpose.

The fourth section is devoted to the problem of existence of finite-dimensional
simple Filippov superalgebras of
type $A(1,0)$ over an algebraically closed field of cha\-racteristic zero. In
the fifth one we treat an analogous problem for the type $A(0,n)$ with $n$ $\in \mathbb{N}$.
In each of these two final sections we restrict our considerations to the case of the Lie
superalgebra that gives the name to the type and solve the existence problem of the mentioned skewsymmetric
homomorphisms. It turns out that the required homomorphisms do not exist.
Therefore, there are no simple finite-dimensional Filippov superalgebras
of type  $A(0,n)$ over an algebraically closed field of characteristic $0$, as stated in the main result of this article (Theorem \ref{princ}). Moreover, as a corollary of its proof, we see that there is no simple finite-dimensional Filippov superalgebra ${\cal F}$ of type $A(0,n)$ such that ${\cal F}$ is a highest weight module over $A(0,n)$.

In what follows, by $\Phi$ we denote  an algebraically closed field of
 characteristic $0$, by $F$ a field of characteristic $0$, by $k$
 a field and by $\left< w_{\upsilon}; \ \upsilon \in \Upsilon \right>$
  a linear space
spanned by the family of vectors
 $\{w_{\upsilon}; \ \upsilon \in \Upsilon \}$ over a field (the field is clear from the context). The symbol $:=$ denotes an
 equality by definition.


\section{Reduction to Lie superalgebras}

From now on,  we denote by $\F$ an $n$-ary Filippov superalgebra.
Let us denote by $\F^*$ ($L(\F)$) the associative (Lie) superalgebra
generated  by the operators $L(x_1,\ldots,x_{n-1})$, $x_i\in \F$.
The algebra $L(\F)$ is called {\it the algebra of multiplications} of $\F$.

\begin{lemma} \textnormal{\cite{Poj8}} \it Let $\F=\F\0\oplus\F\1$ be a simple
finite-dimensional Filippov superalgebra over a field of characteristic $0$
with  $\F\1\neq 0$. Then $L=L(\F)=L\0\oplus L\1$ has nontrivial even and odd
parts.
\end{lemma}

\begin{theorem} \textnormal{\cite{Poj8}} \it If $\F$ is a simple
finite-dimensional Filippov superalgebra over a field of
characteristic $0$, then $L=L(\F)$ is a semisimple Lie superalgebra.
\label{teo1}
\end{theorem}

Given an $n$-ary superalgebra $A$ with a multiplication
$(\cdot,\cdots,\cdot)$, we have $End(A)=End\0A\oplus End\1A$. The element $D\in
End_{\bar s}A$ is called a {\it derivation} of degree $s$ of $A$ if, for every
$a_1,\ldots,a_n\in A, p(a_i)=p_i$, the following equality holds
$$D(a_1,\ldots,a_n)=\sum_{i=1}^{n}(-1)^{sq_i}(a_1,\ldots,Da_i,\ldots,a_n),$$
where $q_i=\sum_{j=1}^{i-1}p_j$.
We denote by $Der_{\bar s}A \subset End_{\bar s}A$
the subspace of all derivations
of degree $s$ and set $Der(A)=Der\0A\oplus Der\1A$.
The subspace $Der(A)\subset End(A)$ is easily seen to be closed under the
bracket

\begin{equation*} [a,b]=ab-(-1)^{deg(a)deg(b)}ba \end{equation*} (known as the \textit{supercommutator}) and it is called
{\it the superalgebra of derivations} of $A$.

\sloppy
Fix $n-1$ elements $x_1,\ldots,x_{n-1}\in A$, $i\in \{1, \ldots, n\}$,
and define a transfor\-mation
$ad_i(x_1,\ldots,x_{n-1})\in End(A)$ by the rule

\begin{equation*} ad_i(x_1,\ldots,x_{n-1})\, x=(-1)^{pq_i}(x_1,\ldots,x_{i-1},x,
x_{i},\ldots,x_{n-1}), \end{equation*} where $p=p(x),p_i=p(x_i),q_i=\sum_{j=i}^{n-1}p_j$.

\sloppy
If, for all $i=1,\ldots,n$ and $x_1,\ldots,x_{n-1}\in A$, the transformations
$ad_i(x_1,\ldots,x_{n-1})\in End(A)$ are derivations of $A$,
then we call them {\it strictly inner derivations} and $A$
 {\it an inner-derivation superalgebra $({\cal ID}$-superalgebra}).
Notice that the $n$-ary Fi\-lippov superalgebras and the $n$-ary commutative
Leibniz algebras are examples of ${\cal ID}$-superalgebras.

Now let us denote by $Inder(A)$ the linear space
spanned by the strictly inner derivations of $A$.
 If $A$ is an $n$-ary $\cal ID$-superalgebra then it is easy to see
 that $Inder(A)$ is an ideal of $Der(A)$.

\begin{lemma} \textnormal{\cite{Poj8}} Given  a simple $\cal ID$-superalgebra $A$ over $k$,
the Lie superalgebra $Inder(A)$ acts faithfully and irreducibly on $A$.
\end{lemma}

Let $\F$ be an $n$-ary Filippov superalgebra over $k$.
Notice that the map $ad:=ad_n:\otimes^{n-1}\F\mapsto Inder(\F)$ satisfies
$$[D,ad(x_1,\ldots,x_{n-1})]=\sum_{i=1}^{n-1}(-1)^{pq_i}
ad(x_1,\ldots,x_{i-1},Dx_i,x_{i+1},\ldots,x_{n-1}),$$
for all $D\in Inder(\F)$, and the associated map
$(x_1,\ldots,x_n)\mapsto ad(x_1,\ldots,x_{n-1})\, x_n$
from $\otimes^{n}\F$ to $\F$ is $\mathbb{Z}_2$-skewsymmetric.
If we regard $\F$ as an $Inder(\F)$-module
then $ad$ induces an $Inder(\F)$-module morphism from
the $(n-1)$-th exterior power $\wedge^{n-1}\F$ to
$Inder(\F)$ (which we also denote by $ad$) such that the map
 $(x_1,\ldots,x_n)\mapsto ad(x_1,\ldots,x_{n-1})\, x_n$ is $\mathbb{Z}_2$-skewsymmetric.
(Note that in $\wedge^{n-1}\F$ we
have$:$ $x_1\wedge\ldots\wedge x_i\wedge x_{i+1}\wedge\ldots \wedge x_{n-1}=
-(-1)^{p_ip_{i+1}}x_1\wedge\ldots\wedge x_{i+1}\wedge x_{i}\wedge\ldots
\wedge x_{n-1}$.)
Conversely, if $(L,V,ad)$ is a triple with $L$ a Lie superalgebra, $V$ an
$L$-module, and $ad$ an $L$-module morphism from $\wedge^{n-1}V \mapsto L$
 such that the map
 $(v_1,\ldots,v_n)\mapsto  ad(v_1\wedge\ldots\wedge v_{n-1})\, v_n$
 from $\otimes^{n}V$ to $V$ is $\mathbb{Z}_2$-skewsymmetric
(we call the homomorphisms of this type {\it skewsymmetric}),
 then $V$ becomes an  $n$-ary Filippov superalgebra by defining
\begin{center} $[v_1,\ldots,v_n]=ad(v_1\wedge\ldots\wedge v_{n-1})\, v_n.$ \end{center}
Therefore, we obtain a correspondence between the set of $n$-ary
 Filippov superalgebras and the set of triples $(L,V,ad)$,
 satisfying the conditions above.

We assume that all vector spaces appearing in the following are
finite-dimensional over $F$.

If $\F$ is a simple $n$-ary Filippov superalgebra then Theorem \ref{teo1}
shows that the Lie superalgebra $Inder(\F)$  is semisimple, and  $\F$
 is a faithful and irreducible $Inder(\F)$-module. Moreover,
 the $Inder(\F)$-module morphism  $ad:\wedge^{n-1}\F\mapsto Inder(\F)$
 is surjective.

Conversely, if $(L,V,ad)$ is a triple such that $L$ is a semisimple
 Lie superalgebra over $F$, $V$ is a
faithful irreducible $L$-module, $ad$ is a surjective
 $L$-module morphism from $\wedge^{n-1}V$ onto the adjoint module $L$,
 and the map
 $(v_1,\ldots,v_n)\mapsto ad(v_1\wedge\ldots\wedge v_{n-1})\, v_n$
 from $\otimes^{n}V$ to $V$ is $\mathbb{Z}_2$-skewsymmetric, then the corresponding
 $n$-ary Filippov superalgebra is simple. A triple with these conditions
 will be called a {\it good triple}. Thus, the problem of determining
 the simple  $n$-ary Filippov superalgebras over $F$ can be
 translated to that of finding the good triples.

\section{Some notations and results on Lie superalgebras}

In this section, we recall some notations and results from \cite{Kac}
on the Lie superalgebra $A(0,n)$ (and its irreducible faithful
finite-dimensional representations). We also give some explicit constructions
which we shall use some later in the study of the simple finite-dimensional Filippov superalgebras of type $A(0,n)$.
Let us start recalling the definition
of induced module.

Let $\L$ be a Lie superalgebra, $U(\L)$ its universal enveloping superalgebra
\cite{Kac}, $H$ a subalgebra of $\L$, and $V$ an $H$-module. The module
$V$ can be extended to $U(H)$-module. We consider the $\mathbb{Z}_2$-graded space
$U(\L)\otimes_{U(H)}V$ (this is the quotient space of $U(\L)\otimes V$ by the
linear span of the elements of the form $gh\otimes v-g\otimes h(v)$,
$g\in U(\L)$, $h\in U(H)$). This space can be endowed with the structure of a
$\L$-module as follows$:$ $g(u\otimes v)=gu\otimes v,g\in \L,u\in U(\L),v\in V$.
The so-constructed $\L$-module is said to be {\it induced from the $H$-module
$V$} and is denoted by $Ind_H^{\L}V$.

From now on, we denote by $G$  a contragredient Lie superalgebra over
$\Phi$ and  consider it with the ``standard'' $\mathbb{Z}$-grading
\cite[Sections 5.2.3 and 2.5.7]{Kac}.

Let $G=\oplus_{i\geq -d}G_i$.
 Set $H=(G_0)\0=\left< h_1,\ldots,h_n\right>$, $N^+=\oplus_{i>0}G_i$ and
 $B=H\oplus N^+$. Let $\Lambda\in H^*, \Lambda(h_i)=a_i\in\Phi$,
$\left<v_\Lambda \right>$ be an one-dimensional $B$-module for which
$N^+(v_\Lambda)=0, h_i(v_\Lambda)=a_iv_\Lambda$. Let $\delta_i \in H^*$, $\delta_i(h_j)=\delta_{ij}$ where
$\delta_{ij}$ is Kronecker's delta.
Let $V_\Lambda=Ind_B^G\left<v_\Lambda \right>/I_{\Lambda}$,
where $I_{\Lambda}$ is the (unique) maximal submodule of the $G$-module
$V_\Lambda$. Then $\Lambda$ is called
the {\it highest weight} of the $G$-module
 $V_\Lambda$.
 By \cite{Kac}, every faithful irreducible finite-dimensional
 $G$-module may be obtained in this manner. Note that the condition
$1\otimes v_\Lambda\in V\0(V\1)$ gives a $\mathbb{Z}_2$-grading on $V_\Lambda$.

\begin{lemma} \textnormal{\cite{Poj8}} \label{soma} \it
 Let $V$ be a module over a Lie superalgebra $G$,
 let $V=\oplus V_{\gamma_i}$ be its weight decomposition, and let $\phi$ be
 a homomorphism from $\wedge^mV$ into $G$. Then, for all $v_i\in
 V_{\gamma_i}$,

\begin{center}
\begin{tabular} {ll}
 $\phi(v_1,\ldots,v_m)\in G_{\gamma_1+\ldots+\gamma_m}$,
& \hspace{0,5cm} if \ $\gamma_1+\ldots+\gamma_m$ is a root of $G$, \\
$\phi(v_1,\ldots,v_m)=0$, & \hspace{0,5cm} otherwise.
\end{tabular}
\end{center}
 \end{lemma}

Let $G$ be a contragredient Lie superalgebra of rank $n$,
$U=Ind_B^G\left<v_\Lambda \right>$, and
$V=V_\Lambda=U/N$
be a finite-dimensional representation of $G$,
where $N=I_{\Lambda}$ is a maximal proper submodule of the $G$-module
$V_\Lambda$.
Let $G=\oplus_{\alpha}G_{\alpha}$
 be a root decomposition of $G$ relative to a Cartan subalgebra $H$.
Denote by ${\cal A}$ the following set of roots$:$
${\cal A}=\{\alpha:g_{\alpha}\notin B\}$.

\begin{lemma} \textnormal{\cite{PS}} \label{pesos} \it Let $g_{\alpha}\in G_{\alpha}$ and $g_{\alpha}\otimes
v\neq 0$ $(v=v_\Lambda)$. Then

\begin{center} $g_{\alpha}^j\otimes v\in U_{\sum_{i=1}^n
(j\alpha(h_i)+\Lambda(h_i))\delta_i}$ \end{center} for all $j\in{\mathbb N}$, and there
exists a minimal
 $ k\in {\mathbb N}$ such that  $g_{\alpha}^k\otimes v \in N$. Moreover, the set
${\cal E}_{\alpha,k}=\{1\otimes v,g_{\alpha}\otimes v,\ldots,
g_{\alpha}^{k-1}\otimes v\}$ is linearly independent in $V$.
Setting $h=[g_{-\alpha},g_{\alpha}]$, we have

\begin{enumerate}
\item $\Lambda(h)=-\frac{(k-1)\alpha(h)}{2}$ if either $g_{\alpha}\in G\0$
   or $k\notin 2{\mathbb N}$;

\item $ \alpha(h)=0$ if $g_{\alpha}\in G\1$ and $k\in 2{\mathbb N}$.
\end{enumerate}
 \end{lemma}

\begin{remark}
Note that if we start with a root $\beta$ then there exists
$s\in{\mathbb N}$ such that ${\cal E}_{\beta,s}$ is linearly independent,
but ${\cal E}_{\alpha,k}\cup {\cal E}_{\beta,s}$
may not be linearly independent.
\end{remark}

Recall that a set ${\cal E}$ is called a {\it pre-basis} of a vector space $W$
if $\left<{\cal E}\right>=W$.

Let $\{g_{\alpha_1}^{k_1}\ldots g_{\alpha_s}^{k_s}\otimes v; k_i\in {\mathbb
N}_0, \alpha_i\in{\cal A}\}$ be a pre-basis of $U$. As we have seen above, for
every $i=1,\ldots,s$, there exists a minimal $p_i\in {\mathbb N}$ such that
$g_{\alpha_i}^{p_i}\otimes v \in N$. Using the induction on the word length, it is easy
to show that $\{g_{\alpha_1}^{k_1}\ldots g_{\alpha_s}^{k_s}\otimes v; k_i\in
{\mathbb N}_0,k_i<p_i,\alpha_i\in{\cal A}\}$ is a pre-basis of $U/N$.

We finish this part with some more notations that we use in the two next sections:

\hspace{-0,7cm} $\bullet$ the symbol $\= $
denotes an equality up to a nonzero coefficient;

\hspace{-0,7cm} $\bullet$ $\underline{u,v}_{\ t}$ means that the elements $u$ and
$v$ are $t$-times repeating $\underbrace{u,v,\ldots,u,v}_{2t}$\,, being
the index $t$ omitted when its value is clear from the context.

\section{Simple Filippov superalgebras of type $A(1,0)$}

Consider the Lie superalgebra $G=A(1,0):=sl(2,1)$. It consists of the matrices of type
$$\left(
\begin{array}{cc|c}
a         &    *         &     *    \\
*         &   b  &     *    \\ \hline
*&  * &     c   \\
\end{array}
\right),
$$
 where $a+b=c$ (\textit{i.e.}, the \textit{supertrace} of these matrices is zero) and all the entries lie in the ground field. We have the following elements in $G
$:
\begin{eqnarray*}
&&\left.
\begin{array}{l}
 h_1 =\ e_{11}-e_{22},\ h_2=e_{22}+e_{33}\\
 e_{12}=g_{\epsilon_1-\epsilon_2} \in G_{\epsilon_1-\epsilon_2},
 e_{21}=g_{\epsilon_2-\epsilon_1} \in G_{\epsilon_2-\epsilon_1},
\end{array}
\right\}\in G\0\\
&&\left.
\begin{array}{l}
e_{13}\in G_{-\epsilon_2}, e_{31}\in G_{\epsilon_2}, \\
e_{32}\in G_{\epsilon_1}, e_{23}\in G_{-\epsilon_1}
\end{array}
\right\}\in G\1.
\end{eqnarray*}

The space $H:= G_0=\left< h_1,h_2\right>$ is a Cartan subalgebra of
$A(1,0)$ and $\epsilon_i,\ i=1,2$, are the linear functions on
$H$ such that $\epsilon_i(e_{jj})=\delta_{ij}$, where
$\delta_{ij}$ is Kronecker's delta. Then
$\Delta=\Delta_0\cup\Delta_1$ is a root system for $A(1,0)$, where
$\Delta_0=\{ 0; \epsilon_1-\epsilon_2;
\epsilon_2-\epsilon_1\}$, and $\Delta_1=\{\pm\epsilon_j;j=1,2\}$.
The roots $\{ \epsilon_1-\epsilon_{2},-\epsilon_1\}$
 are simple.

We have the following
standard grading of $A(1,0)$ \cite[Section 5.2.3]{Kac}:

\begin{equation*} G=\sum_{i=-2}^{2}G_i= \{G_{\epsilon_2}\}
\oplus\{G_{\epsilon_2-\epsilon_1},G_{\epsilon_1}\}\oplus
\{G_0\}\oplus\{G_{\epsilon_1-\epsilon_2},G_{-\epsilon_1}\}\oplus
\{G_{-\epsilon_2}\}. \end{equation*} Because of this, the set
\begin{eqnarray}\label{exp2} {\cal E}=\Big\{\
g_{\epsilon_2}^{t_2}
g_{\epsilon_2-\epsilon_1}^{\gamma}g_{\epsilon_1}^{t_1}
\otimes v\ :\ \gamma \in{\mathbb N}_0, t_i \in \{0,1\}\Big\} \label{expo}
\end{eqnarray}
is a basis of the induced module
$M=Ind_B^G\left<v_\Lambda \right>$ ($v=v_\Lambda$).

Let $V$ be an irreducible module over $G=A(1,0)$ with the highest
weight $\Lambda$, $\Lambda(h_i)=a_i$. Denote $\Lambda$ by $(a_1,a_2)$.
Applying Lemma \ref{pesos}, it is easy to conclude that

\begin{center} $a_1\in {\mathbb N}_0$. \end{center}

If $u\in V_{\gamma}$ (or $G_{\gamma}$) then we may write $h_ku=p_k(u)u$
($[h_k,u]=p_k(u)u$) and call {\it $k$-weight} of $u$ to $p_k(u)$.

In this section, the symbol $w \oti{i}{j}$ means that
$p_1(w \otimes v)=i$ and $p_2(w \otimes v)=j$ (the same with the notation $\mathop u\limits_{j}^{i}$).

\begin{lemma} \label{lema}
Let $V=V_{\Lambda}$ be an irreducible module over $A(1,0)$ with $\Lambda=(1,a_2)$. Then
\begin{center} $g_{\epsilon_2-\epsilon_1}^2\otimes v=g_{\epsilon_1}^2 \otimes v=g_{\epsilon_2}^2 \otimes v=0$. \end{center}
\end{lemma}

\proof As $a_1 \neq 0$ then $g_{\epsilon_2-\epsilon_1} \otimes v \neq 0$. So, by Lemma \ref{pesos},
$g_{\epsilon_2-\epsilon_1}^2 \otimes v = 0$. As far as $g_{\epsilon_1}$ and $g_{\epsilon_2}$, notice that
$t_i \in \{0,1\}$  in (\ref{expo}). In fact, for example, we deduce that $g_{\epsilon_1}^2 \otimes v=0$ from
$[g_{\epsilon_1},g_{\epsilon_1}]=0$. \qed

\begin{corollary} \label{segundabase}
Under the assumptions in Lemma \ref{lema}, the following set is
a pre-basis of $V$$:$
\begin{center} $\{1 \oti{1}{a_2}, g_{\epsilon_2-\epsilon_1} \oti{-1}{1+a_2}, g_{\epsilon_1}\oti{2}{a_2},
g_{\epsilon_2}\oti{0}{1+a_2},
g_{\epsilon_1}g_{\epsilon_2} \oti{1}{1+a_2}, g_{\epsilon_2-\epsilon_1}g_{\epsilon_1} \oti{0}{1+a_2}, g_{\epsilon_2-\epsilon_1}g_{\epsilon_2} \oti{-2}{2+a_2}\}$. \end{center}
 \end{corollary}

\begin{lemma}
Let $V=V_{\Lambda}$ be an irreducible module over $A(1,0)$ with $\Lambda=(0,a_2)$. Then$:$
\begin{center} $g_{\epsilon_2-\epsilon_1}\otimes v=g_{\epsilon_1}^2 \otimes v=g_{\epsilon_2}^2 \otimes v=0$. \end{center}
\end{lemma}

\proof
Suppose that $g_{\epsilon_2-\epsilon_1} \otimes v \neq 0$. By Lemma \ref{pesos}, there exists a minimal $k \in \mathbb{N}$ such that $g_{\epsilon_2-\epsilon_1}^k \otimes v =0$. Concretely,
$k=1$, which is a contradiction. Thus, $g_{\epsilon_2-\epsilon_1} \otimes v = 0$, as we wanted to prove. It only remains to take into account (\ref{expo}).\qed

\begin{corollary} \label{base}
With the assumptions of the previous lemma, the following set is a pre-basis of $V$$:$
\begin{center} $\{1 \oti{0}{a_2}, g_{\epsilon_1}\oti{1}{a_2}, g_{\epsilon_2}\oti{-1}{1+a_2},
g_{\epsilon_1}g_{\epsilon_2} \oti{0}{1+a_2} \}$. \end{center}
 \end{corollary}

\begin{theorem} There exist no simple finite-dimensional Filippov
superalgebras of type $A(1,0)$ over $\Phi$. \label{teoremaA(1,0)}\end{theorem}

\proof Let $G=A(1,0)$, $V$ be a finite-dimensional faithful irreducible module
over $G$ with the highest weight $\Lambda=(a_1, a_2)$ \ ($a_1\neq 0$)
and $\phi$ be a surjective skewsymmetric homomorphism from $\wedge^sV$ on $G$.
Then there exist $u_i\in V_{\gamma_i}$ such that
\begin{eqnarray}
  \label{primeiro}
  \phi(u_1,\ldots,u_{s})&=&g_{\epsilon_2-\epsilon_1}.
\end{eqnarray} By Lemma \ref{soma}, $\sum_{i=1}^sp_1(u_i)=-2$. From Lemma \ref{pesos}, we can conclude that
  $g_{\epsilon_2-\epsilon_1}^{a_1}\otimes v\neq 0$.
Thus, $\phi(u_1,\ldots,u_{s})(g_{\epsilon_2-\epsilon_1}^{a_1-1}\otimes v)
\neq 0$. Since $\phi$
is a skewsymmetric homomorphism,
we have $\phi(u_1,\ldots,u_{i-1},g_{\epsilon_2-\epsilon_1}^{a_1-1}\otimes v,
u_{i+1},\ldots,u_{s})\neq 0$.
As $p_1(g_{\epsilon_2-\epsilon_1}^{a_1-1}\otimes v)=-a_1+2$, the inequality
$|p_1(u_i)+a_1|\leq 2$ follows.  Consequently,
$\phi$ does not exist if $a_1\geq 4$.

Consider \framebox{$a_1=3$}. Thus, by (\ref{primeiro}) and from the fact that $p_1(u_i)<0$, we deduce

\begin{center} $\phi(\stackrel{-1}{u_1},\stackrel{-1}{u_2})=g_{\epsilon_2-\epsilon_1}$. \end{center} From the nonzero action on $1 \otimes v$, we arrive at $\phi(\stackrel{-1}{u_1},1 \oti{3}{})\=
g_{\epsilon_1-\epsilon_2}$. Acting on $g_{\epsilon_2-\epsilon_1}\otimes v$, we
have $\phi(g_{\epsilon_2-\epsilon_1}\oti{1}{ },1\oti{3}{ })\neq 0$,
which is a weight contradiction.

Now let us take \framebox{$a_1=2$}. According to (\ref{primeiro}) and taking into account that $p_1(u_i) \in \{-4,-3,-2,-1,0\}$, we only have two possibilities$:$
\begin{center} i) $\phi(\stackrel{0}{u_1}, \ldots, \stackrel{0}{u_{s-1}}, \stackrel{-2}{u_s})
=g_{\epsilon_2-\epsilon_1}$; \ \ \ \ \ \ \ \ \ \ \ \ \ ii) $\phi(\stackrel{0}{u_1}, \ldots, \stackrel{0}{u_{s-2}},\stackrel{-1}{u_{s-1}}, \stackrel{-1}{u_s})
=g_{\epsilon_2-\epsilon_1}$. \end{center}
Consider i). Let us suppose first that $1 \otimes v$ is even.
Acting on $1 \otimes v$,
we have $\phi(
\ox{u_1}{0}{},\ldots,\ox{u_{s-1}}{0}{},\ox{1 \otimes v}{2}{})\=
g_{\epsilon_1-\epsilon_2}$.
Since the action of $g_{\epsilon_1-\epsilon_2}$ on $g_{\epsilon_2-\epsilon_1} \otimes v$ provides a
non zero element, we arrive at
$\phi(g_{\epsilon_2-\epsilon_1}
\ox{\otimes}{0}{}v,\ox{u_{2}}{0}{},\ldots,\ox{u_{s-1}}{0}{},\ox{1 \otimes v}{2}{})\=
g_{\epsilon_1-\epsilon_2}$. Finally, acting once again on $g_{\epsilon_2-\epsilon_1} \otimes v$,
we obtain $\phi(g_{\epsilon_2-\epsilon_1} \otimes v, u_3,\ldots, u_{s-1}, g_{\epsilon_2-\epsilon_1} \otimes v) \neq 0$ and, from here, a skewsymmetric contradiction.
As far as i), it only remains to assume that $1 \otimes v$  is odd. As before,
acting on $1 \otimes v$,
we have $\phi(
\ox{u_1}{0}{},\ldots,\ox{u_{s-1}}{0}{},\ox{1 \otimes v}{2}{})\=
g_{\epsilon_1-\epsilon_2}$. Then, acting repeatedly on $g_{\epsilon_2-\epsilon_1} \otimes v$,
we get $\phi(\underline{g_{\epsilon_2-\epsilon_1} \otimes v})\neq 0$. Since
$p_1(g_{\epsilon_2-\epsilon_1} \otimes v)=0$, we have $\phi(\underline{g_{\epsilon_2-\epsilon_1} \otimes v}) \=\alpha h_1 + \beta h_2$, with nonzero right-hand side. From the multiplication by
 $g_{-\epsilon_1}$ we get $0=\alpha$. On the other hand, the multiplication by
 $g_{-\epsilon_2}$ leads to $\alpha=\beta$. Henceforth, $\alpha=\beta=0$ which is a contradiction. In the subcase
 ii), the multiplication
by $g_{\epsilon_1-\epsilon_2}$
gives either
\begin{center} $\phi(
\ox{u_1}{0}{},\ldots,\ox{v_{i}}{2}{},\ldots,\ox{u_{s-2}}{0}{},
\ox{u_{s-1}}{-1}{},\ox{u_s}{-1}{})(1\oti{2}{})\neq 0$, for some $v_i$ with $i \in \{1,\ldots,s-2\}$, \\
or $\phi(
\ox{u_1}{0}{},\ldots,\ox{u_{s-2}}{0}{},
\ox{v_{s-1}}{1}{},\ox{u_s}{-1}{})(1\oti{2}{})\neq 0$, for some $v_{s-1}$. \end{center}
In both cases, replacing $u_s$ by $1 \otimes v$, we arrive at a
weight contradiction.

Consider \framebox{$a_1=1$}.

1) Assume that $a_2 \neq 0$. We study two subcases, 1.1)
and 1.2), in what follows.

1.1) Suppose that $a_2 \neq -1$. Assume also that
$1 \otimes v$ is even and $\phi(u_1,\ldots,u_s)=g_{\epsilon_2-\epsilon_1}$. Taking into
account the weight considerations and Corollary \ref{segundabase}, we can suppose that
\begin{equation} \phi(g_{\epsilon_2-\epsilon_1} \oti{-1}{}, g_{\epsilon_2-\epsilon_1} \oti{-1}{},
\mathop\protect{\overline{u_3,\ldots,u_s}}\limits^{0})=g_{\epsilon_2-\epsilon_1} \label{poss1} \end{equation}

\begin{equation} \textnormal{or} \
\phi(g_{\epsilon_2-\epsilon_1}g_{\epsilon_2}\oti{-2}{}, \mathop\protect{\overline{u_2,\ldots,u_s}}\limits^{0})=
g_{\epsilon_2-\epsilon_1}. \label{poss2} \end{equation}
In (\ref{poss1}), as $g_{\epsilon_2-\epsilon_1} \otimes v$ is even, we have arrived to a skewsymmetry
contradiction. As far as (\ref{poss2}), acting on $1 \otimes v$, we have
\begin{center}
$\phi(1 \oti{1}{},\mathop\protect{\overline{u_2,\ldots,u_s}}\limits^{0}) \= g_{\epsilon_1} (g_{-\epsilon_2})$.
\end{center}
If  $\phi(1 \oti{1}{},\mathop\protect{\overline{u_2,\ldots,u_s}}\limits^{0}) \= g_{\epsilon_1}$
 then, acting once again on $1 \otimes v$ (nonzero action since $a_2 \neq 0$), we arrive at
 a skewsymmetry contradiction. Suppose now that we have
 $\phi(1 \oti{1}{},\mathop\protect{\overline{u_2,\ldots,u_s}}\limits^{0}) \=g_{-\epsilon_2}$.
 If all elements $u_2,\ldots,u_s$ have $1-$weight equal to $0$ then their $2-$weight is equal to $1+a_2$.
 Thus, $a_2+(s-1)(1+a_2)=-1$, that is, $a_2=-1$,
  a contradiction. So, we can assume that one of the elements $u_2,\ldots,u_s$
  has $1-$weight less than zero ($-1$ or $-2$); without loss of generality, suppose that it is $u_2$. Then, acting on $g_{\epsilon_2} \otimes v$ and substituting the mentioned element, we arrive at

\begin{center} $\phi(1 \oti{1}{}, g_{\epsilon_2} \oti{0}{}, \mathop\protect{\overline{u_3,\ldots,u_s}}\limits^{\textnormal{$1$}}) \neq 0$,
\end{center} (we considered $p_1(u_2)=-1$ since the other possibility leads to an immediate weight contradiction) where one of the elements in $\{u_3,\ldots,u_s\}$ has $1-$weight equal to $1$. If this element is $1 \otimes v$ then we obtain a skewsymmetry contradiction. If not, then we have

\begin{center} $\phi(1 \oti{1}{}, g_{\epsilon_2} \oti{0}{}, g_{\epsilon_1}g_{\epsilon_2} \oti{1}{},
\mathop\protect{\overline{u_4,\ldots,u_s}}\limits^{0})
\= g_{\epsilon_1-\epsilon_2}$. \end{center} If all the elements $u_4,\ldots,u_s$ have $1-$weight equal to $0$ then
they all have $2-$weight equal to $1+a_2$. So, as above, $a_2=-1$, which is a contradiction. Therefore,
we may assume that $p_1(u_4) \in \{-2, -1\}$. From the action on $g_{\epsilon_2} \otimes v$ a weight
contradiction arises. Consider $1 \otimes v$ odd. Assume that $\phi(u_1,\ldots,u_s)=
g_{-\epsilon_1}$. Once again by Corollary \ref{segundabase}, suppose first that

\begin{center} $\phi(g_{\epsilon_2-\epsilon_1} \oti{-1}{},\mathop\protect{\overline{u_2,\ldots,u_s}}\limits^{0})
=g_{-\epsilon_1}$. \end{center} Acting
on $g_{\epsilon_1} \otimes v$, we arrive at $\phi(g_{\epsilon_1} \oti{2}{},
\mathop\protect{\overline{u_2,\ldots,u_s}}\limits^{0}) \= g_{\epsilon_1-\epsilon_2}$.
Then, from the nonzero action on $g_{\epsilon_2} \otimes v$,
we can write $\phi(g_{\epsilon_1} \oti{2}{},g_{\epsilon_2} \oti{0}{},
\mathop\protect{\overline{u_3,\ldots,u_s}}\limits^{\textnormal{$0$ or $1$ or $2$}}) \neq 0$.
If we have a weight contradiction then we are done; if not, we obtain $
\phi(g_{\epsilon_1} \oti{2}{},g_{\epsilon_2} \oti{0}{},
\mathop\protect{\overline{u_3,\ldots,u_s}}\limits^{\textnormal{$0$}}) \= g_{\epsilon_1-\epsilon_2}$,
and we just have to act on $g_{\epsilon_2} \otimes v$ one more time to get a skewsymmetry contradiction. Assume now that

\begin{center} $\phi(u_1,g_{\epsilon_2-\epsilon_1}g_{\epsilon_2} \oti{-2}{},
\mathop\protect{\overline{u_3,\ldots,u_s}}\limits^{0})=g_{-\epsilon_1}$,
where $u_1 \in \{1\oti{1}{}, g_{\epsilon_1}g_{\epsilon_2} \oti{1}{}\}$. \end{center}
Then we can act on $g_{\epsilon_1} \oti{2}{}$ to obtain a weight contradiction.

1.2) Suppose that $a_2 = -1$ and $\phi(u_1,\ldots,u_s)=g_{-\epsilon_1}$.
Then, by Corollary \ref{segundabase}, we have to analise three situations. Assume first that

\begin{center} $\phi(g_{\epsilon_2-\epsilon_1} \oti{-1}{},
\mathop\protect{\overline{u_2,\ldots,u_s}}\limits^{0})=g_{-\epsilon_1}$.
\end{center} From the nonzero action on $g_{\epsilon_1} \otimes v$, we obtain
$\phi(g_{\epsilon_1} \oti{2}{}, \mathop\protect{\overline{u_2,\ldots,u_s}}\limits^{0}) \= g_{\epsilon_1-\epsilon_2}$.
If one of the elements $u_2,\ldots,u_s$ has negative $1-$weight then, acting on $g_{\epsilon_2} \oti{0}{}$,
we arrive at a weight contradiction. Otherwise, we have
$\phi(g_{\epsilon_1} \oti{2}{-1}, \mathop {u_2}\limits_{0}^{0},
\ldots, \mathop {u_s}\limits_{0}^{0}) \= g_{\epsilon_1-\epsilon_2}$. Acting on $g_{\epsilon_2-\epsilon_1} \otimes v$,
we obtain
$\phi(g_{\epsilon_1} \oti{2}{-1}, g_{\epsilon_2-\epsilon_1} \oti{-1}{0},
\mathop {u_3}\limits_{0}^{0}, \ldots, \mathop {u_s}\limits_{0}^{0}) \=
g_{-\epsilon_2}$. Finally, by the nonzero action on $g_{\epsilon_1}g_{\epsilon_2}\oti{1}{}$, we obtain a weight contradiction. Suppose now that

\begin{center} $\phi(u_1,g_{\epsilon_2-\epsilon_1}g_{\epsilon_2} \oti{-2}{},
\mathop\protect{\overline{u_3,\ldots,u_s}}\limits^{0})=g_{-\epsilon_1}$, where $u_1 \in \{1 \oti{1}{},
g_{\epsilon_1}g_{\epsilon_2} \oti{1}{}\}$. \end{center} Then, acting on $g_{\epsilon_1} \oti{2}{}$,
we get a weight contradiction.

2) Assume that $a_2 =0$. Thus, by Corollary \ref{segundabase},
the $2-$weights of the elements that constitute the pre-basis of $V$ are greater or equal to zero.
Hence, it is impossible to have $\phi(u_1,\ldots,u_s) = g_{\epsilon_1-\epsilon_2}$, for some
$u_1, \ldots, u_s$.

Suppose now that \framebox{$a_1=0$ and $a_2 \neq 0$}.
Assume that $\phi(u_1,\ldots,u_s)=g_{\epsilon_2-\epsilon_1}$. By Corollary
\ref{base}, we may write: $u_1=u_2=g_{\epsilon_2} \otimes v$. Therefore, $1 \otimes v$ is even.
Note that $a_2 \neq -1$ since otherwise we cannot find $g_{\epsilon_2-\epsilon_1}$ in the image of $\phi$.
Suppose that  $\phi(v_1,\ldots,v_s)=\alpha h_1 + \beta h_2$, where $\beta \neq 0$.
Then $v_i \notin \{1 \otimes v, g_{\epsilon_1}g_{\epsilon_2} \otimes v\}$. Denote
$\underline{g_{\epsilon_1} \otimes v,g_{\epsilon_2} \otimes v}_k$ by $x_k$, where $k \in \mathbb{N}$.
Thus, we may assume that

\begin{center} $\phi(x_k)=\alpha h_1 + \beta h_2, \ \beta \neq 0$. \end{center} Taking into account the $2$-weights,
we conclude that $a_2 = -1/2$. We have $g_{\epsilon_1} \phi(x_k)=-\alpha g_{\epsilon_1}$ and
$g_{\epsilon_2} \phi(x_k)=(\alpha-\beta) g_{\epsilon_2}$. Thus,

\begin{equation} -k\phi(g_{\epsilon_1} \otimes v, g_{\epsilon_1}g_{\epsilon_2} \otimes v,x_{k-1})=-\alpha g_{\epsilon_1},
\label{um} \end{equation}

\begin{equation} k \phi(g_{\epsilon_2} \otimes v, g_{\epsilon_1}g_{\epsilon_2} \otimes v, x_{k-1})=(\alpha-\beta) g_{\epsilon_2}. \label{dois} \end{equation}
From (\ref{um}) and (\ref{dois}), acting on $g_{\epsilon_2} \otimes v$ and $g_{\epsilon_1} \otimes v$, respectively,
we arrive at

\begin{center} $k[g_{\epsilon_1}g_{\epsilon_2} \otimes v,x_k]=-\alpha g_{\epsilon_1}g_{\epsilon_2} \otimes v$
 and $-k[g_{\epsilon_1}g_{\epsilon_2} \otimes v,x_k]=(\beta-\alpha)g_{\epsilon_1}g_{\epsilon_2} \otimes v$. \end{center}
Whence, $\beta=2\alpha$. On the other hand, acting in (\ref{um}) on $g_{\epsilon_2} \otimes v$,
we get $k \phi(x_k) g_{\epsilon_1}g_{\epsilon_2} \otimes v=-\alpha g_{\epsilon_1}g_{\epsilon_2} \otimes v$, \textit{i.e.},

\begin{center} $k(\alpha h_1 + \beta h_2) g_{\epsilon_1}g_{\epsilon_2} \otimes v=-\alpha g_{\epsilon_1}g_{\epsilon_2}
\otimes v$. \end{center} Hence $k \beta=-2\alpha$. Therefore $\alpha=\beta=0$, a contradiction.

To finish the proof, consider \framebox{$a_1 = a_2 = 0$}.
By Corollary \ref{base}, $\{1 \oti{}{0}, g_{\epsilon_1}\oti{}{0}, g_{\epsilon_2}\oti{}{1},
g_{\epsilon_1}g_{\epsilon_2} \oti{}{1} \}$ is a pre-basis of $V$.
As $p_2(g_{\epsilon_1-\epsilon_2})=-1$, it is impossible to find $u_1, \ldots, u_s \in V$
such that $\phi(u_1,\ldots,u_s)=g_{\epsilon_1-\epsilon_2}$. \qed

\begin{corollary} \label{corol}
There are no simple finite-dimensional Filippov
superalgebras of type $A(0,1)$ over $\Phi$.
\end{corollary}

\proof It suffices to take into account that $A(m,n) \simeq A(n,m)$, \cite[Section 4.2.2]{Kac}.
\qed

\section{Simple Filippov superalgebras of type $A(0,n)$}

Recall that $A(m,n):=sl(m+1,n+1)$ for $m \neq n$ and $m, n \in \mathbb{N}_0$. This Lie superalgebra consists of the matrices of type
$$\left(
\begin{array}{c|c}
A         &   B \\
\hline
&\\[-3mm]
C    &  D\\
\end{array}
\right),
$$
 where $A \in M_{(m+1)\times(m+1)}(F), B \in M_{(m+1)\times (n+1)}(F), C \in M_{(n+1) \times (m+1)}(F), D \in M
 _{(n+1)\times (n+1)}(F)$ and $tr(A)=tr(D)$. Let us write some elements in $G=A(m,n)$:
\begin{eqnarray*}
&&\left.
\begin{array}{l}
 h_i =\ e_{ii}-e_{i+1,i+1}, \ \ i=1,\ldots,m,m+2,\ldots,m+n+1,\\
 h_{m+1}=e_{m+1,m+1}+e_{m+2,m+2},\\
 e_{kl}:=g_{\epsilon_k-\epsilon_l} \in G_{\epsilon_k-\epsilon_l}, \ \ k,l=1,\ldots,m+1 \ \textnormal{or} \ k,l=m+2,\ldots,m+n+2,
\end{array}
\right\}\in G\0\\
&&\left.
\begin{array}{l}
e_{kl}:=g_{\epsilon_k-\epsilon_l} \in G_{\epsilon_k-\epsilon_l}, \ \ k=1,\ldots,m+1, l=m+2,\ldots,m+n+2, \\
\hspace{3,6cm} \textnormal{or} \
k=m+2,\ldots,m+n+2, l=1,\ldots,m+1.
\end{array}
\right\}\in G\1
\end{eqnarray*}

The space $H:=G_0=\left< h_1,\ldots,h_{m+n+1}\right>$
is a Cartan subalgebra
of $A(m,n)$, and $\epsilon_i$
are the linear functions
on $H$ defined by its values on $h_1,\ldots,h_{m+n+1}$
and the conditions $\epsilon_i(e_{jj})=\delta_{ij}$, where $\delta_{ij}$
is Kronecker's delta. Then $\Delta=\Delta_0\cup\Delta_1$ is
a root system for $A(m,n)$, where $\Delta_0=\{ 0; \epsilon_k-\epsilon_l, k,l=1,\ldots,m+1 \ \textnormal{or} \ k,l=
m+2,\ldots,m+n+2\}$,
and $\Delta_1=\{\epsilon_k-\epsilon_l, k=1,\ldots,m+1, l=m+2, \ldots,m+n+2 \ \textnormal{or} \ k=m+2,\ldots,m+n+2, l=1,\ldots,m+1\}$.
The roots $\{\alpha_i:= \epsilon_i-\epsilon_{i+1},i=1,\ldots,m+n+1\}$
 are simple.

The conditions $deg \ g_{\alpha_i}=1, deg \ g_{-\alpha_i}=-1$ give us the standard grading of $A(m,n)$, \cite[Section 5.2.3]{Kac}. The negative part of this grading is $G_{\epsilon_k-\epsilon_l}, l < k$. Because of this, the set
\begin{eqnarray}\label{exp}{\cal E}=\left\{
\prod_{l < k} g_{\epsilon_k-\epsilon_l}^{\gamma_{kl}}
\otimes v\ :\ \gamma_{kl} \in{\mathbb N}_0\right\}
\end{eqnarray}
is a pre-basis of the induced module
$M=Ind_B^G\left<v_\Lambda \right>$ (where $v=v_\Lambda$).

Let $V$ be an irreducible module over $G=A(m,n)$ with the highest
weight $\Lambda$, $\Lambda(h_i)=a_i$. Denote $\Lambda$ by $(a_1,\ldots,a_{m+n+1})$.
Applying Lemma \ref{pesos}, we have

\begin{equation} a_i \in \mathbb{N}_0 \ \ \textnormal{if} \ i \neq m+1. \label{am+1}\end{equation}

In what follows, we assume that $m=0$ and, because of Corollary \ref{corol}, $n \in \mathbb{N} \setminus \{1\}$.

\begin{lemma} \label{jueves}
Let $V=V_{\Lambda}$ be an irreducible module over $G=A(0,n)$ with the highest weight $\Lambda=
(a_1,\ldots,a_{n+1})$, where $a_1>0$ and $\sum_{i=2}^{n+1} a_i= 1$.
Assume that $(G,V,\phi)$ is a good triple. Then $\frac{1}{2} \leq a_1 \leq 1$.
\end{lemma}

\proof
Suppose that $\phi(u_1,\ldots,u_s)=g_{\epsilon_2-\epsilon_1}$. By the action on $1 \otimes v$,
we obtain $p_1(u_i) \geq -1+a_1$. If $a_1 > 1$ then $p_1(u_i) > 0$. Thus, we have a contradiction since
$p_1(g_{\epsilon_2-\epsilon_1})=0$.
On the other hand, through the same action and with $h=e_{11}+e_{n+2,n+2}$, we have $p_h(u_i) \leq -1 +a_1$. If
$a_1 < \frac 1 2$ then $p_h(u_i) < -\frac 1 2$, which contradicts $p_h(g_{\epsilon_2-\epsilon_1})=-1$.
\qed

\begin{lemma} \label{hope}
Let $V=V_\Lambda$ be an irreducible module over $A(0,2)$ with the highest weight
$\Lambda=(1,1,0)$. Then there are no good triples $(A(0,2),V,\phi)$.
\end{lemma}

\proof
Assume that $\phi(u_1,\ldots,u_s)=g_{\epsilon_2-\epsilon_1}$, $H=e_{22}-e_{44}$, and $h=e_{11}+e_{44}$.
 By the action on $1 \otimes v$, we have
$0\leq p_2(u_i), p_H(u_i) \leq 4$ and $-2\leq p_h(u_i)\leq 0$ with $i=1,2$.
If $s\geq 3$ then we may assume that $p_2(u_3)=p_H(u_3)=0$ and $u_3\mapsto 1 \otimes v$
gives a weight contradiction.

Note that $g_{\epsilon_3-\epsilon_1}\otimes v\neq 0,$ since otherwise $g_{\epsilon_1-\epsilon_2}
g_{\epsilon_3-\epsilon_1}\otimes v= 0,$ and $g_{\epsilon_3-\epsilon_2}\otimes v= 0.$

Suppose that $\phi(u_1,u_2)=g_{\epsilon_3-\epsilon_2}$. We also have
$g_{\epsilon_3-\epsilon_2}
g_{\epsilon_2-\epsilon_1}\otimes v\neq 0.$
Therefore, the action on $g_{\epsilon_2-\epsilon_1}\otimes v$ gives
$p_h(u_i)=0,i=1,2.$ If $p_2(u_1)\leq -2$ then $u_1\mapsto g_{\epsilon_2-\epsilon_1} \otimes v$
gives a $(h_2,h)$-contradiction (a contradiction considering the $2$-weights and the $h$-weights of the elements).
If $p_2(u_1)=p_2(u_2)=-1$ then $u_1\mapsto g_{\epsilon_2-\epsilon_1} \otimes v$ leads to
$u=\phi(g_{\epsilon_2-\epsilon_1} \otimes v,u_2)\in \{g_{\epsilon_2-\epsilon_1},g_{\epsilon_2-\epsilon_4}\}$,
considering $h_2, h$-weights. If $u\=g_{\epsilon_2-\epsilon_1}$ then $u_2\mapsto 1 \otimes v$
gives a $h_2$-contradiction. If $u\=g_{\epsilon_2-\epsilon_4}$ then $u_2\mapsto g_{\epsilon_4-\epsilon_2} \otimes v$
allows us to arrive at $\phi(g_{\epsilon_2-\epsilon_1} \otimes v,g_{\epsilon_4-\epsilon_2} \otimes v)
\=g_{\epsilon_2-\epsilon_3}$
and $g_{\epsilon_4-\epsilon_2} \otimes v\mapsto g_{\epsilon_3-\epsilon_1} \otimes v$ leads to a
$h$-contradiction. \qed

\begin{lemma} \label{mi}  Let $V=V_{\Lambda}$ be an irreducible module over $G=A(0,n)$ with the highest weight
$\Lambda=(a_1, 0, \ldots, 0)$. Assume that $(G,V,\phi)$ is a good triple. Then $0 < a_1 \leq 1$.
\end{lemma}

\proof Let $\beta=\epsilon_i-\epsilon_j$ with $i > j$ and $j \neq 1$.
Thus, $g_\beta \in G_{\0}$.
By Lemma \ref{pesos}, $g_\beta \otimes v =0$. Consider
$g_{\epsilon_i-\epsilon_1}$ with $i \neq 1$.
Notice that $g_{\epsilon_i-\epsilon_1} \in G_{\1}$. Then, from
$[g_{\epsilon_i-\epsilon_1}, g_{\epsilon_i-\epsilon_1}]\otimes v=0$, we conclude that $g_{\epsilon_i-\epsilon_1}^2
\otimes v =0$. Henceforth,

\begin{equation} \label{last}
\mathcal{B}=\left\{\prod_{i=2}^{n+2} g_{\epsilon_i-\epsilon_1}^{\gamma_i} \otimes v: \gamma_i \in \{0,1\}\right\}
\end{equation} is a pre-basis of $V$.

Observe that if $a_1 \leq 0$ then the $1$-weights of the elements in $\mathcal{B}$
are negative or zero. Thus, it is impossible to have $\phi(u_1,\ldots,u_s)=
g_{\epsilon_1-\epsilon_3}$.

Now suppose that $a_1 >1$.
Let $\phi(u_1,\ldots,u_s)=g_{\epsilon_{n+2}-\epsilon_1}$ and $h=e_{11}+e_{n+2,n+2}$.
 The action on $1 \otimes v$ leads to $|a_1-p_h(u_i)| \leq 1$ and $p_h(u_i)>0$ with $i=1, \ldots, s$, a contradiction.
\qed

 \begin{theorem} \label{princ} There are no simple finite-dimensional Filippov
superalgebras of type $A(0,n)$ over $\Phi$. \end{theorem}

\proof
Suppose that $V$ is a finite-dimensional irreducible module
over $G=A(0,n)$ with the highest weight $\Lambda=(a_1, \ldots, a_{n+1})$ ($a_2+\ldots+a_{n+1}\neq 0$),
and
 $\phi$ is a surjective skewsymmetric homomorphism from $\wedge^sV$ on $G$. Then there exist $u_i\in V_{\gamma_i}$ such that
\begin{eqnarray}
  \label{eqnova2}
  \phi(u_1,\ldots,u_{s})&=&g_{\epsilon_{n+2}-\epsilon_2}.
\end{eqnarray} Let $H=h_2+\ldots+h_{n+1}=e_{22}-e_{n+2,n+2}$. By Lemma \ref{soma}, $\sum_{i=1}^sp_H(u_i)=-2$. From Lemma \ref{pesos},
  $g_{\epsilon_{n+2}-\epsilon_2}^{a_2+\ldots+a_{n+1}}\otimes v\neq 0$.
 Since $\phi$
is a skewsymmetric homomorphism,
$\phi(u_1,\ldots,u_{i-1},g_{\epsilon_{n+2}-\epsilon_2}^{a_2+\ldots+a_{n+1}-1}\otimes v,
u_{i+1},\ldots,u_{s})\neq 0$.
As $p_H(g_{\epsilon_{n+2}-\epsilon_2}^{a_2+ \ldots+a_{n+1}-1}\otimes v)=2-a_2-\ldots-a_{n+1}$, the inequality
$|p_H(u_i)+a_2+\ldots+a_{n+1}|\leq 2$ follows.  From here we see that
the required skewsymmetric homomorphism does not exist if $a_2+\ldots+a_{n+1}\geq 4$.

Consider the case \framebox{$a_2+\ldots+a_{n+1}=3$}. In this case, $p_H(u_i) <0$.
So, by (\ref{eqnova2}), we have $\phi(\stackrel{-1}{u_1},\stackrel{-1}{u_2})=g_{\epsilon_{n+2}-\epsilon_2}$ and, acting on $1 \otimes v$, we arrive at $\phi(\stackrel{-1}{u_1}, 1\oti{3}{})\=
g_{\epsilon_{2}-\epsilon_{n+2}}$. Acting on $g_{\epsilon_{n+2}-\epsilon_2}\otimes v$, we
obtain $\phi(g_{\epsilon_{n+2}-\epsilon_2}\oti{1}{ },1\oti{3}{ })\neq 0$,
which is a weight contradiction.

Now let us take \framebox{$a_2+\ldots+a_{n+1}=2$}. As $\sum_{i=1}^s p_H(u_i)=-2$ and, in this case, $p_H(u_i) \leq 0$, we can only have

\begin{center} i) $\phi(\stackrel{-2}{u_1}, \stackrel{0}{u_{2}}, \ldots,  \stackrel{0}{u_s})
=g_{\epsilon_{n+2}-\epsilon_2}$ \ \ \ or \ \ ii) $\phi(\stackrel{-1}{u_1}, \stackrel{-1}{u_2},
\stackrel{0}{u_3}, \ldots, \stackrel{0}{u_s})
=g_{\epsilon_{n+2}-\epsilon_2}$. \end{center}
First consider i). Let us suppose that $1 \otimes v$ is even.
Acting on $1 \otimes v$,
we have $\phi(1 \oti{2}{},
\ox{u_2}{0}{},\ldots,\ox{u_{s}}{0}{})\=
g_{\epsilon_2-\epsilon_{n+2}}$.
Then, acting twice on $g_{\epsilon_{n+2}-\epsilon_2} \otimes v$, we arrive at
$\phi(\underline{g_{\epsilon_{n+2}-\epsilon_2} \otimes v}_{2}, u_3, \ldots, u_s) \neq 0$ which leads to a skewsymmetry
contradiction.
Assume now that $1 \otimes v$  is odd. Then, acting on $1 \otimes v$ and, repeatedly, on $g_{\epsilon_{n+2}-\epsilon_2}\otimes v$, we get
$\phi(\stackrel{2}{1 \otimes v},\underline{g_{\epsilon_{n+2}-\epsilon_2}\stackrel{0}{\otimes}v})
\= g_{\epsilon_2-\epsilon_{n+2}}$. From here, analizing the $1$-weights, we
conclude that $a_1 = 1$. Let

\begin{equation} \phi(u_1,\ldots,u_s)=g_{\epsilon_2-\epsilon_1}. \label{Litta}\end{equation} By the action on $1 \otimes v$, we obtain $|3-p_H(u_i)| \leq 2$. Thus, (\ref{Litta}) can not occur.
Consider ii). The multiplication
by $g_{\epsilon_2-\epsilon_{n+2}}$
gives either $$\phi(\stackrel{-1}{u_1}, \stackrel{-1}{u_2},
\stackrel{0}{u_3}, \ldots, \stackrel{2}{w},\ldots, \stackrel{0}{u_s})
(1\oti{2}{})\neq 0\  \mbox{ or }\
\phi(\stackrel{-1}{u_1}, \stackrel{1}{v_2},
\stackrel{0}{u_3}, \ldots, \stackrel{0}{u_s})(1\oti{2}{})\neq 0, \ \textnormal{for some $w, v_2$}.$$
In both cases, replacing $u_1$ by $1\otimes v$, we arrive at a
weight contradiction.

Consider the case \framebox{$a_2+\ldots+a_{n+1}=1$}. Let $h=e_{11}+e_{n+2,n+2}$.

I) Suppose that $a_1 <0$ and

\begin{equation}  \phi(u_1,\ldots,u_s)=g_{\epsilon_{n+2}-\epsilon_1}. \label{desp} \end{equation} Through the action on $1 \otimes v$, we get $|-p_h(u_i)+a_1-1|\leq 1$. Thus, $p_h(u_i) <0$, being (\ref{desp}) impossible. So, $a_1 \geq 0$.

II) Assume that $a_1>0$. From Lemma \ref{jueves}, we have $\frac 1 2 \leq a_1 \leq 1$.

IIa) Consider $a_1 \neq 1$.
 Let $\phi(u_1,\ldots,u_s)=g_{\epsilon_{2}-\epsilon_{1}}$. By the action on $1
\otimes v$, we have $p_{H}(u_i) \in \{0,1,2,3,4\}$. Consequently, after the action on $1 \otimes v$, we arrive at

\begin{center} $\phi(\stackrel{1}{u_1}, 1 \oti{1}{},
\stackrel{0}{u_3},\ldots,\stackrel{0}{u_s}) \= g_{\epsilon_{2}-\epsilon_{n+2}}$.\end{center} Replacing every $u_k (k \geq
3)$ by $g_{\epsilon_{n+2}-\epsilon_1} \otimes v$ and acting one more time on the mentioned element, we get
$\phi(1 \oti{1}{}, g_{\epsilon_{n+2}-\epsilon_1} \oti{0}{}, \ldots, g_{\epsilon_{n+2}-\epsilon_1} \oti{0}{}) \neq 0$.
Analyzing the $H$, $h$ and $2$-weights involved, we conclude that $a_2=0$ and

\begin{center} $\phi(1 \oti{1}{}, g_{\epsilon_{n+2}-\epsilon_1} \oti{0}{}, \ldots, g_{\epsilon_{n+2}-\epsilon_1}
\oti{0}{}) \= g_{\epsilon_i-\epsilon_{n+2}}, \ i\neq 1, 2, 3, n+2$. \end{center}
Through the multiplication by $g_{\epsilon_{n+2}-\epsilon_1}$, we obtain $\phi(\underline{g_{\epsilon_{n+2}-\epsilon_1}
\otimes v}) \= g_{\epsilon_{i}-\epsilon_1}, i \neq 1, 2, 3, n+2$.
Trough the consecutive analises of the
$3, \ldots, n$-weights, we eliminate all the possibilities for $i$.

IIb) Take $a_1=1$. 1) Suppose that $a_2=1$. Assume that $\phi(u_1, \ldots, u_s)=g_{\epsilon_2-\epsilon_1}$. From the nonzero
action on $1 \otimes v$, we conclude that $p_H(u_i), p_2(u_i) \in \{0,\ldots,4\}, p_1(u_i) \in \{0,1,2\}$ and
$p_h(u_i) \in \{-2,-1,0\}$. Putting the $h$-weights above the elements and the $1$-weights underneath
them, we have $\phi(\mathop {u_1}\limits_{0}^{-1}, \mathop {u_2}\limits_{0}^{0},
\ldots, \mathop {u_s}\limits_{0}^{0})=g_{\epsilon_2-\epsilon_1}$. The mentioned action
implies

\begin{equation} \phi(1 \oti{0}{1}, \mathop {u_2}\limits_{0}^{0},
\ldots, \mathop {u_s}\limits_{0}^{0}) \= g_{\epsilon_1-\epsilon_{n+2}} (g_{\epsilon_2-\epsilon_i}, i \neq 1,2, n+2)
. \label{Quinta} \end{equation} Consider the first possibility in (\ref{Quinta}). Observe that
 $p_2(g_{\epsilon_2-\epsilon_1})=1, p_2(1 \otimes v)=1, p_2(g_{\epsilon_1-\epsilon_{n+2}})=0$. Thus,
$p_2(u_1)=2$ and $\sum_{j=2}^s p_2(u_j)=-1$, a contradiction. Assume the second possibility. Notice that
$p_2(g_{\epsilon_2-\epsilon_i})$ is either $1$ (when $i \neq 3$) or $2$ (if $i=3$). If $i \neq 3$ then $p_2(u_1)=1$ and
we arrive at the contradiction
$\phi(\mathop{\mathop{u_1}\limits_{0}^{-1}}\limits_{1},
\mathop{1 \oti{0}{1}}\limits_{1}, \mathop{\mathop{u_3}\limits_{0}^0}\limits_{0},\ldots,
\mathop{\mathop{u_s}\limits_{0}^0}\limits_{0}) \neq 0$,
where the weight lines refer to $h, h_1, h_2$, respectively. Take $i=3$. On one hand, as $p_H(g_{\epsilon_2-\epsilon_1})=
1, p_H(1 \otimes v)=1$ and $p_H(g_{\epsilon_2-\epsilon_3})=1$, we have $p_H(u_1)=1$.
On the other hand, assuming $h'=h_2+h_3=e_{22}-e_{44}$, as
$p_{h'}(g_{\epsilon_2-\epsilon_1})=
1, p_{h'}(1 \otimes v)=1$ and $p_{h'}(g_{\epsilon_2-\epsilon_3})=1$, we have $p_{h'}(u_1)=1$. Thus, we
 get the weight contradiction
$\phi(\mathop {u_1}\limits_{1}^{1}, 1 \oti{1}{1}, \mathop {u_3}\limits_{0}^{0},
\ldots, \mathop {u_s}\limits_{0}^{0}) \neq 0$, where the $H$-weights are above the elements and the $h'$-weights are
underneath them, for $n \neq 2$. By Lemma \ref{hope}, $n=2$ can not occur.

2) Suppose that $a_{n+1}=1$. In this subcase, we put the $(n+1)$-weights above the elements and the
$1$-weights underneath them.
 Assume that $\phi(u_1, \ldots, u_s)=g_{\epsilon_2-\epsilon_1}$. From the action on
$1 \otimes v$, we conclude that $p_1(u_i) \in \{0,1,2\}, p_{n+1}(u_i) \in \{-1,0,1,2,3\}, p_h(u_i) \in \{-2,-1,0\}$.
If there is a $j$ such that $p_{n+1}(u_j) <0$ then we arrive at a weight contradiction obtained
from $u_j \mapsto 1 \oti{1}{1}$. So, $p_{n+1}(u_i) \geq 0$. Therefore, the nonzero action on
 $1 \otimes v$ leads to
$\phi(\mathop{\mathop{u_1}\limits_0^0}\limits_{-1},1 \mathop{\oti{1}{1}}\limits_{\hspace{-0,2cm} 0},
\mathop{\mathop{u_3}\limits_0^0}\limits_0,\ldots,\mathop{\mathop{u_s}\limits_0^0}\limits_0) \neq 0$, where the
third weight line refers to $h$. Thus, we have

\begin{equation} \label{Covi1}
\phi(\mathop{\mathop{u_1}\limits_0^0}\limits_{-1},1 \mathop{\oti{1}{1}}\limits_{\hspace{-0,2cm} 0},
\mathop{\mathop{u_3}\limits_0^0}\limits_0,\ldots,\mathop{\mathop{u_s}\limits_0^0}\limits_0) \=
g_{\epsilon_2-\epsilon_{n+2}}.
\end{equation} From (\ref{Covi1}), through the action on $g_{\epsilon_{n+2}-\epsilon_2} \otimes v$,
we arrive at the weight contradiction
$\phi(g_{\epsilon_{n+2}-\epsilon_2}\mathop{\oti{0}{0}}\limits_{\hspace{-0,3cm}1},1 \mathop{\oti{1}{1}}
\limits_{\hspace{-0,2cm} 0},
\mathop{\mathop{u_3}\limits_0^0}\limits_0,\ldots,\mathop{\mathop{u_s}\limits_0^0}\limits_0) \neq 0$.

3) Suppose that there is a $2 <j <n+1$ such that $a_j=1$. Here we put the $1$-weights above the
elements and the $j$-weights underneath them. Assume that $\phi(u_1,\ldots,u_s)=
g_{\epsilon_{j+1}-\epsilon_j}$. From the action on $1 \otimes v$, we conclude that $p_j \in \{-3,-2,-1,0,1\},
p_1(u_i) \in \{0,1,2\}$ and $p_{j+1} \in \{-1,0,1,2,3\}$. Suppose that $p_j(u_1)=1$. Then,
through the mentioned action, we get $\phi(1 \oti{1}{1},\mathop{\underline{\stackrel{0}{u_2},
\ldots, \stackrel{0}{u_s}}}\limits_{-3}) \neq 0$.
As this is a weight contradiction then $p_j(u_i) \leq 0$. Consider

\begin{equation*} \phi(\mathop{\stackrel{0}{u_1}}\limits_{0},\ldots,\mathop{\stackrel{0}{u_{s-2}}}\limits_0
\mathop{\stackrel{0}{u_{s-1}}}\limits_{-1}, \mathop{\stackrel{0}{u_s}}\limits_{-1})=g_{\epsilon_{j+1}-\epsilon_j}.
\end{equation*}  Multiplying by $g_{\epsilon_j-\epsilon_{j+1}}$ and acting on $1 \otimes v$, we arrive at
$\phi(\mathop\protect{\underline{\stackrel{0}{u_1},\ldots,\stackrel{0}{u_{s-1}}}}\limits_{1},1 \oti{1}{1}) \neq 0$,
one more
weight contradiction. Now take

\begin{equation} \label{Covi2} u=\phi(\mathop{\stackrel{0}{u_1}}\limits_{0},\ldots,
\mathop{\stackrel{0}{u_{s-1}}}\limits_{0}, \mathop{\stackrel{0}{u_s}}\limits_{-2})=g_{\epsilon_{j+1}-\epsilon_j}.
\end{equation} The action on $1 \otimes v$ leads to $w=\phi(\mathop{\stackrel{0}{u_1}}\limits_{0},\ldots,
\mathop{\stackrel{0}{u_{s-1}}}\limits_{0},1 \oti{1}{1}) \neq 0$. Thus, $w \in \{g_{\epsilon_1-\epsilon_{j+1}},
 g_{\epsilon_2-\epsilon_{j+1}}\}$. Notice that $p_{j+1}(u)=1, p_{j+1}(1 \otimes v)=0$ and $p_{j+1}(w)=-1$.
 Consequently, $p_{j+1}(u_s)=2$ and we may assume that $p_{j+1}(u_1)=-1$. So, from (\ref{Covi2}), acting on $1 \otimes v$, we get the weight contradiction
 $\phi(\mathop{1 \oti{1}{1}}\limits_{0},
\mathop{\underline{\mathop {u_2}\limits_{0}^{0},\ldots,\mathop {u_{s-1}}\limits_{0}^{0}}}\limits_{0},
\mathop{\mathop{u_s}\limits_{-2}^0}\limits_{2}) \neq 0$, where the third weight line refers to $h_{j+1}$.

III) Assume that $a_1=0$. Consider $h$ as previously defined. Suppose that
$\phi(u_1, \ldots, u_s)=
g_{\epsilon_{n+2}-\epsilon_1}$. The action on $1 \otimes v$ allows us to conclude that $p_1(u_i),
p_h(u_i) \in \{-2,-1,0\}$. Thus, we have

\begin{equation}
\phi(\mathop {u_1}\limits_{0}^{-1},\mathop {u_2}\limits_{0}^{0},\ldots,\mathop {u_s}\limits_{0}^{0})=
g_{\epsilon_{n+2}-\epsilon_1}, \label{Covi3} \end{equation} where, here and in what follows,
we consider the $1$-weights above the elements and the $h$-weights underneath them.
By the action on $1 \otimes v$ in (\ref{Covi3}), we arrive at

\begin{equation} \label{Covi4} \phi(\mathop {u_1}\limits_{0}^{-1}, 1 \oti{0}{-1},
\mathop {u_3}\limits_{0}^{0},\ldots,\mathop {u_s}\limits_{0}^{0}) \=
g_{\epsilon_{j}-\epsilon_1}, \ j \neq 1, 2, n+2. \end{equation} If $g_{\epsilon_j-\epsilon_1} \otimes v \neq 0$
then we get the weight contradiction
$\phi(1 \oti{}{-1},1 \oti{}{-1},\mathop {u_3}\limits_{0}^{},\ldots, \mathop {u_s}\limits_{0}^{}) \neq 0$. So,
$g_{\epsilon_j-\epsilon_1} \otimes v = 0$. Then there exists $k > j$ such that $g_{\epsilon_j-\epsilon_1}
g_{\epsilon_k-\epsilon_j} \otimes v \neq 0$.
If $k \neq n+2$ then, through the action on $g_{\epsilon_k-\epsilon_j} \oti{}{-1}$, we obtain a weight contradiction.
Thus, we may assume that $a_{n+1}=1$. Replace all $u_t$ ($t \geq 3$) in (\ref{Covi4})
by $g_{\epsilon_{n+2}-\epsilon_{j_r}} \oti{0}{0}$ and act one more time on this element. We obtain

\begin{equation} u:=\phi(1 \oti{0}{-1}, g_{\epsilon_{n+2}-\epsilon_{j_1}} \oti{0}{0},
g_{\epsilon_{n+2}-\epsilon_{j_3}} \oti{0}{0}, \ldots,
g_{\epsilon_{n+2}-\epsilon_{j_s}} \oti{0}{0})\= g_{\epsilon_2-\epsilon_1} (g_{\epsilon_i-\epsilon_{n+2}}),
\end{equation} with  $i \neq 1, 2, n+2$. If $u \=g_{\epsilon_2-\epsilon_1}$ then we arrive at the
contradiction $0 \= g_{\epsilon_2-\epsilon_{n+2}}$ through the multiplication by $g_{\epsilon_1-
\epsilon_{n+2}}$. In the latter occasion, notice that $j_r < n+1$. In fact, if there is a $q=j_t=n+1$, for some $t$,
then the analysis
of the weights over $e_{11}+e_{qq}$ leads to a weight contradiction. Moreover, taking $i= n+1, n, \ldots$,
the consecutive analysis of the $i$-weights implies weight contradictions that allow us to exclude all the
possibilities for $i$.

Consider the case \framebox{$a_2+\ldots+a_{n+1}=0$}. So, $a_2=\ldots=a_{n+1}=0$. Let $h=e_{11}+e_{n+2,n+2}$.
From
Lemma \ref{mi}, we have $0 < a_1 \leq 1$.

I) Take $a_1=1$. In what follows, let us put the $h$-weights and the $1$-weights above and underneath the
elements, respectively. Assume that $\phi(u_1, \ldots, u_s)=g_{\epsilon_2-\epsilon_1}$. The action on $1 \otimes v$ allows us to conclude that $p_h(u_i) \in \{-1,0,1\}$ and $p_1(u_i) \in \{0,1,2\}$. Thus, we have

\begin{center}
$\phi(\mathop {u_1}\limits_0^{-1},\mathop\protect{\overline{\mathop{u_2}\limits_{0},\ldots,\mathop{u_s}\limits_{0}}}
\limits^{\textnormal{$0$}}) = g_{\epsilon_2-\epsilon_1}$.
\end{center} Suppose that $p_h(u_2)=-1$. Then, through the action on $1 \otimes v$, we obtain

\begin{center} $\phi(\mathop {u_1}\limits_0^{-1},1 \oti{1}{1},
\mathop\protect{\overline{\mathop{u_3}\limits_{0},\ldots,\mathop{u_s}\limits_{0}}}
\limits^{1}) \= g_{\epsilon_1-\epsilon_i}$, $\quad i \neq 1, 2, n+2$. \end{center} The action on $g_{\epsilon_i
-\epsilon_1} \otimes v$ leads to the weight contradiction $\phi(g_{\epsilon_i
-\epsilon_1} \oti{0}{},1 \oti{1}{},
\mathop\protect{\overline{\mathop{u_3},\ldots,\mathop{u_s}}}
\limits^{1}) \neq 0$. So, consider now

\begin{center}
$\phi(\mathop {u_1}\limits_0^{-1},\mathop{u_2}\limits_{0}^0,\ldots,\mathop{u_s}\limits_{0}^0) = g_{\epsilon_2-\epsilon_1}$.
\end{center} By the action on $1 \otimes v$, we arrive at

\begin{center} $\phi(1 \oti{1}{1},\mathop{u_2}\limits_{0}^0,\ldots,\mathop{u_s}\limits_{0}^0) \= g_{\epsilon_1-\epsilon_t}$, $\quad t \neq 1, 2, n+2$. \end{center} From the repeated action on $g_{\epsilon_{t_r-\epsilon_1}} \oti{0}{0}$, we get $\phi(1 \otimes v,
g_{\epsilon_{t_2-\epsilon_1}} \otimes v, \ldots, g_{\epsilon_{t_s-\epsilon_1}}
\otimes v) \= g_{\epsilon_1-\epsilon_i}$, with $i \neq 1, 2, n+2$. The weights over $e_{11}+e_{t_2t_2}$ give
a weight contradiction.

II) Suppose now that $0 < a_1 < 1$. Consider

\begin{equation} \label{aqui} \phi(u_1, \ldots, u_s)=g_{\epsilon_3-\epsilon_1}. \end{equation}

IIa) Assume first that $n \neq 2$. Let $h=e_{11}+e_{n+2,n+2}$.
Notice that $g_{\epsilon_4-\epsilon_1} g_{\epsilon_2-\epsilon_1} \otimes v \neq 0$. From
(\ref{aqui}) and the nonzero action on $g_{\epsilon_4-\epsilon_1} g_{\epsilon_2-\epsilon_1} \otimes v$,
we have $|-3-p_h(u_i)+a_1| \leq 1$. As $p_h(u_i)=\tau_i + a_1$, where $\tau_i \in \mathbb{Z}$, then
$|-3-\tau_i| \leq 1$. Taking into account that $p_h(g_{\epsilon_3-\epsilon_1})=-1$ and $0 < a_1 <1$, it is impossible to have
(\ref{aqui}).

IIb) Assume now that $n=2$. Let $h=e_{11}+e_{44}$ and $h_3=e_{33}-e_{44}$. Once again from
(\ref{aqui}) and the nonzero action on $g_{\epsilon_4-\epsilon_1} g_{\epsilon_2-\epsilon_1} \otimes v$, we
conclude that $|-2+a_1-p_h(u_i)| \leq 1$. This condition,
as $p_h(g_{\epsilon_3-\epsilon_1})=-1$, leads to $p_h(u_i)=-1+a_1$ for all $i$. Thus,

\begin{center}
$\phi(\mathop {u_1}\limits_k^{-1+a_1}, \mathop{\underline{\stackrel{-1+a_1}{u_2},
\ldots, \stackrel{-1+a_1}{u_s}}}\limits_{1-k})=g_{\epsilon_3-\epsilon_1}$, \quad  $k > 0$,
\end{center}  where the weights over $h$ and $h_3$ are above and underneath the elements, respectively.
The action on $g_{\epsilon_4-\epsilon_1} g_{\epsilon_2-\epsilon_1} \otimes v$ allows us to write
$\phi(g_{\epsilon_4-\epsilon_1}g_{\epsilon_2-\epsilon_1}\oti{-1+a_1}{-1}, \mathop{\underline{\stackrel{-1+a_1}{u_2},
\ldots, \stackrel{-1+a_1}{u_s}}}\limits_{1-k})\neq 0$.
But there is no root element in $A(0,2)$ with both weights over $h$ and $h_3$ being negative. Hence, we obtained
a contradiction. The theorem is proved.
\qed

\begin{corollary}
There is no simple finite-dimensional Filippov superalgebra ${\cal F}$ of type $A(0,n)$ such that ${\cal F}$ is a highest weight module over $A(0,n)$.
\end{corollary}

\renewcommand{\refname}{\begin{center} References \end{center}}

\vspace{0,8cm}

\begin{footnotesize}P. D. Beites \\ Departamento de Matem\'{a}tica and Centro de Matemсtica, Universidade da Beira Interior \\ Covilhу, Portugal
 \\ \textit{E-mail adress}: pbeites@ubi.pt\end{footnotesize}

\vspace{0,4cm}

\begin{footnotesize}A. P. Pozhidaev \\ Sobolev Institute of Mathematics and Novosibirsk State University \\
Novosibirsk, Russia \\ \textit{E-mail adress}: app@math.nsc.ru \end{footnotesize}
\vspace{.5cm}

\end{document}